\documentclass[11pt,preprintnumbers,showpacs,preprint,showpacs]{article}
\usepackage{amsmath}
\usepackage{amssymb}
\usepackage{braket}
\usepackage{graphicx}
\usepackage{dcolumn}
\usepackage{bm}

\title{\bf{A Substitution to Bernoulli Numbers in easier computation of \(\zeta(2k)\) }}
\author{\bf{Srinivasan.A}$^{\dag}$ }

\begin{document}
\maketitle
\begin{abstract}
\noindent  An alternative formula is presented for the evaluation of the zeta function values $\zeta(2k)$ without the need for Bernoulli numbers. Our formula is recursive, and improves the efficiency with which we can calculate large values of the zeta function. 
\end{abstract}

\raggedright
\let\thefootnote\relax\footnotetext {$\dag$ Srinivasan.A (srinivasan1390[at]gmail[dot]com) is a final year undergraduate Engineering student at the National Institute of Technology (NIT), Karnataka(India).}


\section{Introduction}
The Riemann zeta function is highly explored function in the field of mathematics since it lies in the confluence of many fields like number theory, probability, applied statistics, etc. The proof of the Riemann conjecture has impact in the field of security due to it being directly relevant in the distribution of primes. Riemann Conjecture simply stated is as follows \textbf{'All non-trivial zeros of the Zeta function have real part equal to one-half'}.  In this paper we have looked at the already existing formulation\cite{bernoulli1}\cite{riemann1} of the zeta function at even integers. We have aimed in this paper to replace the relation below containing Bernoulli numbers.\\

\begin{flalign*}
\zeta(2k) &=\frac{(-1)^{k+1}B_{2k}(2\pi)^{2k}}{2(2k)!}
\end{flalign*}

\section{Motivation}
	The formula presented above is a simple expressions but is used in the evaluation of zeta function at higher arguments both for the even and the odd case (already existing formulae are shown in the compilation of \cite{bernoulli5}, \cite{bernoulli6}, \cite{bernoulli7} . Although the relation seems straightforward we would like to highlight that computation of Bernoulli numbers for large arguments of zeta function might cause a computational problem. The papers \cite{bernoulli3},\cite{bernoulli4}  contain more efficient algorithms for evaluating Bernoulli numbers for higher values of the argument. Hence instead of employing better algorithms to calculate the Bernoulli numbers faster and accurately we realized that an alternative to the Bernoulli numbers in the computation of $\zeta(2k)$ would be better. Few technical aspects with regard to Bernoulli numbers have been highlighted below.

\section{Bernoulli Numbers}
	Bernoulli numbers are employed in many functions the more prominent ones being the computational aspects of Fermat's Last Theorem, Riemann zeta function, Taylor series expansions of tangent and hyperbolic tangent functions, etc. The general iterative relation used to calculate the bernoulli numbers is presented below.

\raggedright 

\[\frac{1}{r!}\sum_{m=0}^{r-1} \left(
\begin{array}{c}
n \\
r
\end{array}
\right)B_m =  \left\{
\begin{array}{l l}
1 & \quad \mbox{if r=0}  \\

0 & \quad \mbox{if r \textgreater 1}\\

\end{array} \right. \] \\[0.5 cm]

The main problem lies in the computation of higher values of Bernoulli numbers using the above recursive function. Many algorithms for the computation of Bernoulli numbers came up after 1995 when powerful computers and software like Maple, Mathematica failed to compute huge values of Bernoulli numbers. It was reached a conclusion that computation of B(1000)( which has approximately 2000 digits in the decimal representation) would be impossible and hence many researchers started coming up with algorithms for faster and better approximation and calculation of Bernoulli numbers. Although a program was created in 2002 for evaluation of B(4,000,000), this has come up after 6 years of research work. We haven exactly looked into the computational aspects in terms of complexity, time, precision of the formula presented as shown in history. With the simple idea of replacing expression for Bernoulli numbers in the evaluation of $\zeta(2k) $, we have obtained the following result. The crux of the result was the use of Fourier Series to get a result from which we would be able to compare the infinite summation series for zeta function and the Fourier series representation to derive an alternate representation for $\zeta(2k)$ which yields excatly same results.

\section{Main Result}
The general fourier expansion of any function g(x) with the period T.
\begin{flalign}
g(x)= A_o + \sum_{n=0}^{\infty} ( A_n\cos{\frac{2n\pi x}{T}}+ B_n\sin{\frac{2n\pi x}{T}})
\end{flalign}
We consider our as 
\begin{flalign*}
g_k(x) &=(\pi-x)^{2k} \hspace{2pt} s.t. \hspace{2pt} 0 \leq x \leq \pi
\end{flalign*}

\begin{flalign} \notag
A_o &=\frac{1}{\pi} \int_{0}^{\pi} (x-\pi)^{2k} \mathrm{d}x \\ \notag
 &=\frac{1}{\pi}\int_{-\pi}^{\pi} (t)^{2k} \mathrm{d}t \\ 
&=\frac{\pi^{2k}}{2k+1}
\end{flalign}

\begin{flalign*}
A_{n,2k}&=\frac{2}{\pi}\int_{0}^{\pi} (x-\pi)^{2k} \cos{(nx)}  \mathrm{d}x \\[0.2 cm]
&=\frac{2}{\pi}[(x-\pi)^{2k} \frac{\sin{(nx)}}{n} \mid_0^\pi- 2k\int_{0}^{\pi}\frac{sin{(nx)}}{n}(x-\pi)^{2k-1}\mathrm{d}x]\\[0.2 cm]
&=\frac{2}{\pi}[-\frac{2k}{n}[ (x-\pi)^{2k-1} \frac{(-\cos{(nx))}}{n} 
\mid_0^\pi\int_{0}^{\pi}\frac{(-\cos{(nx))}}{n}(2k-1)(x-\pi)^{2k-2}\mathrm{d}x]]\\[0.2cm]
&=\frac{2}{\pi}[-\frac{2k}{n}[ \frac{(-\pi)^{2k-1}}{n}+\frac{2k-1}{n} \int_{0}^{\pi} (\cos{(nx))}(x-\pi)^{2k-2}\mathrm{d}x]]\\[0.2 cm]
&=\frac{2}{\pi}[\frac{2k}{n^2}(\pi)^{2k-1}-\frac{(2k)(2k-1)}{n^2}A_{n,2k-2}\frac{\pi}{2} ]\\[0.2 cm]
&=\frac{4\pi^{2k-2}k}{n^2}-\frac{(2k)(2k-1)}{n^2}A_{n,2k-2}
\end{flalign*}\\[1 cm]

I would like to keep \(a_k=\frac{4k}{n^2\pi^2}\) and \(b_k=\frac{(2k)(2k-1)}{n^2}\)
\begin{subequations}
\begin{flalign} \notag
A_{n,2k}&=\frac{4\pi^{2k-2}k}{n^2}-\frac{(2k)(2k-1)}{n^2}A_{n,2k-2}\\ 
A_{n,2k}&=a_k\pi^{2k}+b_kA_{n,2k-2}\\\notag
A_{n,2k-2}&=a_{k-1}\pi^{2k-2}+b_{k-1}A_{n,2k-4}\\
b_{k}A_{n,2k-2}&=b_{k}a_{k-1}\pi^{2k-2}+b_{k}b_{k-1}A_{n,2k-4}\\\notag
b_{k-1}A_{n,2k-4}&=b_{k-1}a_{k-2}\pi^{2k-4}+b_{k-1}b_{k-2}A_{n,2k-6}\\
b_{k}b_{k-1}A_{n,2k-4}&=b_{k}b_{k-1}a_{k-2}\pi^{2k-4}+b_{k}b_{k-1}b_{k-2}A_{n,2k-6}\\\notag
. . . . . &=. . . . . \\\notag
. . . . . &=. . . . . \\
\prod_{i=0}^{n-1}b_{k-i}A_{2k-2n}&=a_{k-n}\prod_{i=0}^{n-1}b_{k-i}\pi^{2k-2n}+\prod_{i=0}^{n}b_{k-i}A_{n,2k-2n-2}
\end{flalign}
\end{subequations}
Summing the subequations: 
\begin{flalign}
A_{n,2k}&=\sum_{n=0}^{k-2}a_{k-n}\prod_{i=0}^{n-1} b_{k-i}\pi^{2k-2n} +\prod_{i=0}^{k-2}b_{k-i}A_{n,2}
\end{flalign}\\

Using Equation(10) from \textbf{Result 1} with (j=k-2) in equation (4)  

\begin{multline}
A_{n,2k}=\frac{4k}{n^2}\pi^{2k-2}+\frac{4(k-1)}{n^4}\pi^{2k-4}(-2k)(2k-1)+...\\
+\frac{(-2)^{k-1}}{n^{2(k-1)}}(k(k-1)...2(2k-1)(2k-3)...3)\frac{4}{n^2\pi^2} \notag
\end{multline}

\begin{multline}
=\frac{4}{n^2\pi^2}(k\pi^{2k}+\frac{-2(k)(k-1)\pi^{2k-2}}{n^2}+\frac{(-2)^2k(k-1)(k-2)(2k-1)(2k-3)\pi^{2k-4}}{(n^2)^2}+...)\\
+\frac{4}{n^2\pi^2}( (\frac{-2}{n^2})^{k-1}(k(k-1)(k-2)....(2)(1))((2k-1)(2k-3)(2k-5).....(3)(1)))
 \end{multline}\\
From \textbf{Result 2} we get the following generalization.
\begin{flalign}\notag
A_{n,2k}&=\sum_{j=1}^{k}(\frac{4}{n^2\pi^2}(\frac{2(-1)^{j-1}}{n^{2j-2}}k(k-j+1)\frac{(2k-1)!}{(2k-2k+2)!})\pi^{2k-2j+2})\pi^{2(k-(j-1))}\\\notag
&=\frac{8k(2k-1)!\pi^{2k-2}}{n^2\pi^2}\sum_{j=1}^{k}(\frac{(-1)^{j-1}(k-j+1)}{n^{2j-2}\pi^{2j}(2k-2j+2)!})\\
&=\frac{8k(2k-1)!\pi^{2k-2}}{n^2\pi^2}\sum_{j=1}^{k}(\frac{(-1)^{j-1}(k-(j-1))}{n^{2(j-1)}\pi^{2j}(2k-2(j-1)!}))
\end{flalign}
We replace all the (j-1)'s in the above summation by j's and thus summation changes from \(j=1\rightarrow k\) to  \(j=0\rightarrow k-1\)
\begin{flalign}\notag
A_{n,2k}&=\frac{8k(2k-1)!\pi^{2k-2}}{n^2}\sum_{j=0}^{k-1}(\frac{-1}{n^2\pi^2})^j\frac{k-j}{(2k-2j)!}\\\notag
&=\frac{4k(2k-1)!\pi^{2k-2}}{n^2}\sum_{j=0}^{k-1}(\frac{-1}{n^2\pi^2})^j\frac{2(k-j)}{(2k-2j)!}\\
&=\frac{2(2k)!\pi^{2k-2}}{n^2}\sum_{j=0}^{k-1}(\frac{-1}{n^2\pi^2})^j\frac{1}{(2k-2j-1)!}
\end{flalign}\newpage

Substituting (7) and (2) into (1) we get the following. We realize that \(B_{n,2k}\)=0 since \(g_k(x)\) is an odd function.\\

\begin{flalign}\notag
(x-\pi)^{2k}&=\frac{\pi^{2k}}{2k+1}+\sum_{n=1}^{\infty}\frac{2(2k)!\pi^{2k-2}}{n^2}\sum_{j=0}^{k-1}(\frac{-1}{n^2\pi^2})^j\frac{1}{(2k-2j-1)!}cos(nx)\\\notag
With\hspace{10pt} x\hspace{2pt}=\hspace{2pt}0\\\notag
(\pi)^{2k}&=\frac{\pi^{2k}}{2k+1}+\sum_{n=1}^{\infty}\frac{2(2k)!\pi^{2k-2}}{n^2}\sum_{j=0}^{k-1}(\frac{-1}{n^2\pi^2})^j\frac{1}{(2k-2j-1)!}\\\notag
\frac{k}{2k+1}&=\frac{(2k)!\pi^{-2}}{n^2}\sum_{n=1}^{\infty}\sum_{j=0}^{k-1}(\frac{-1}{n^2\pi^2})^j\frac{1}{(2k-2j-1)!}\\
\frac{k}{(2k+1)!}&=\frac{1}{\pi^2}\sum_{j=0}^{k-1}(\frac{-1}{\pi^2})^j\frac{1}{(2k-2j-1)!}\sum_{n=1}^{\infty}(\frac{1}{n^2})^{j+1}
\end{flalign} 
For notational convenience we set
\(\sum_{n=1}^{\infty}(\frac{1}{n^2})^{j}=\zeta(2j)\)
\begin{flalign}\notag
\frac{k}{(2k+1)!}&=\frac{1}{\pi^2}\sum_{j=0}^{k-1}(\frac{-1}{\pi^2})^j\frac{1}{(2k-2j-1)!}\sum_{n=1}^{\infty}(\frac{1}{n^2})^{j+1}\\\notag
&=\frac{1}{\pi^2}(\sum_{j=0}^{k-2}(\frac{-1}{\pi^2})^j\frac{1}{(2k-2j-1)!}\zeta(2j+2)+(\frac{-1}{\pi^2})^{k-1}\zeta(2k))\\\notag
\zeta(2k)(\frac{-1}{\pi^2})^{k-1}+\frac{k}{(2k+1)!}&=-(\sum_{j=0}^{k-2}(\frac{-1}{\pi^2})^{j+1}\frac{1}{(2k-2j-1)!}\zeta(2j+2))\\\notag
\zeta(2k)&=-(\sum_{j=0}^{k-2}(\frac{-1}{\pi^2})^{j+1}\frac{1}{(2k-2j-1)!}\zeta(2j+2)+\frac{k}{(2k+1)!})\pi^{2k}(-1)^k\\
\zeta(2k)=\sum_{n=1}^{\infty}\frac{1}{n^{2k}}&=-(\sum_{j=0}^{k-2}(\frac{-1}{\pi^2})^{j+1}\frac{1}{(2k-2j-1)!}\zeta(2j+2)+\frac{k}{(2k+1)!})\pi^{2k}(-1)^k
\end{flalign}\newpage

 \emph{\bf{Result 1:}}
We try to simplify the product or summation in \[\prod_{i=0}^{j} b_{k-i}\] by finding an equivalent expression for it.
\begin{flalign}\notag
\prod_{i=0}^{j} b_{k-i}&=\\
b_kb_{k-1}b_{k-2}b_{k-3}....b_{k-j}&=\frac{(-2k)(2k-1)}{n^2}\frac{(-2k+2))(2k-3)}{n^2}\frac{(-2k+4)(2k-5)}{n^2}... 
\end{flalign}

\begin{multline}
=\frac{(-1)^{j+1}}{n^{2(j+1)}}2^{j+1}(k(k-1)(k-2)(k-3)....(k-(j-2))(k-(j-1))(k-j))\\
(2k-1)(2k-3)....(2k-2j+1)
\end{multline}\\[1.5 cm]

\emph{\bf{Result 2:}}Let us consider the last term of the  Equation [5]\\[0.5 cm]
\begin{flalign}
\frac{4}{n^2\pi^2}((\frac{-2}{n^2})^{k-1}(k(k-1)(k-2)....(2)(1))(2k-1)(2k-3)&(2k-5).....(3)(1)))\notag\\
(k(k-1)(k-2)....(2)(1))&=k! \notag\\
((2k-1)(2k-3)(2k-5).....(3)(1)) \times \frac{(2k)(2k-2).....(2)}{(2k)(2k-2).....(2)}&=\frac{(2k)!}{2^{k-1}(k!)}\notag
\end{flalign}\\[0.5 cm]

We would like to generalize the Equation [5] by finding the \(j^{th}\) term so as to refer it by an index. Hence the  \(j^{th}\) term of Equation [5] is.\\[0.5 cm]
\begin{flalign*}
\frac{4}{n^2\pi^2}((\frac{-2}{n^2})^{j-1}(k(k-1)....(k-j+1))((2k-1)(2k-3).....(2k-2j-3)))\pi^{2(k-(j-1))}\\
=\frac{4}{n^2\pi^2}(\frac{(-1)^{j-1}2k(k-j+1)}{n^{2(j-1)}}\frac{(2k-1)!}{(2k-2k+2)!})\pi^{2k-2j+2}
\end{flalign*}\\[2 cm]

\newpage
\centering \textbf{Conclusions} \\
\raggedright
Having checked the above result for few values of $\zeta$(2k) using the above formula and verified with \cite{riemann1}\\
\[\zeta(2k)=-(\sum_{j=0}^{k-2}(\frac{-1}{\pi^2})^{j+1}\frac{1}{(2k-2j-1)!}\zeta(2j+2)+\frac{k}{(2k+1)!})\pi^{2k}(-1)^k\]
\begin{flalign*}
k&=1 \rightarrow \zeta(2)=\pi^2/6  \\
k&=2\rightarrow \zeta(4)=-((\frac{-1}{\pi^2})^1\frac{1}{3!}\zeta(2)+\frac{2}{5!})\pi^4=\frac{\pi^4}{90}\\
k&=4 \rightarrow  \zeta(6)=\frac{691\times 2^{11}}{15!}\pi^{12}
\end{flalign*}
The above recursive formula gives as an alternative for the regular Bernoulli dependent formula for the \(\zeta(2k)\) which yields perfect results for all values of $k$ and has avoided the use of Bernoulli numbers.

\raggedright

\end{document}